\documentclass[a4paper,11pt]{amsart}
\usepackage{graphicx}
\usepackage{mathptmx}
\usepackage{amsmath}
\usepackage{amssymb}
\usepackage{enumitem}
\usepackage{xcolor}
\usepackage{xparse}
\NewDocumentCommand{\eulerian}{omm}
 {%
  \genfrac<>{0pt}{}{#2}{#3}%
  \IfValueT{#1}{_{\!#1}}%
 }

 \newmuskip\pFqmuskip

\newcommand*\pFq[6][8]{%
  \begingroup 
  \pFqmuskip=#1mu\relax
  \mathchardef\normalcomma=\mathcode`,
  \mathcode`\,=\string"8000
  \begingroup\lccode`\~=`\,
  \lowercase{\endgroup\let~}\pFqcomma
  {}_{#2}F_{#3}{\left(\genfrac..{0pt}{}{#4}{#5}\bigg|#6\right)}%
  \endgroup
}
\newcommand{\pFqcomma}{{\normalcomma}\mskip\pFqmuskip}

\newtheorem{theorem}{Theorem}
\newtheorem{lemma}[theorem]{Lemma}

\begin{document}

\title[Some identities and properties on degenerate Stirling numbers]{Some identities and properties on degenerate Stirling numbers}

\author{Taekyun  Kim}
\address{Department of Mathematics, Kwangwoon University, Seoul 139-701, Republic of Korea}
\email{tkkim@kw.ac.kr}

\author{DAE SAN KIM}
\address{Department of Mathematics, Sogang University, Seoul 121-742, Republic of Korea}
\email{dskim@sogang.ac.kr}


\maketitle

\begin{abstract}
The aim of this paper is by using generating functions to further study some identities and properties on the degenerate Stirling numbers of the second kind, the degenerate $r$-Stirling numbers of the second kind, the degenerate Stirling numbers of the first kind and the unsigned degenerate Stirling numbers of the first kind.
\end{abstract}

\section{Introduction}
Carlitz [4] obtained some interesting arithmetical and combinatorial results on the degenerate Bernoulli and Euler polynomials and numbers, which are degenerate versions of the Bernoulli and Euler polynomials and numbers.
In recent years, studying degenerate versions of some special numbers and polynomials have drawn the attention of many mathematicians with their regained interests not only in combinatorial and arithmetical properties but also in applications to differential equations, identities of symmetry and probability theory (see [7,9,10,13-16,18,20] and the references therein). These degenerate versions include the degenerate Stirling numbers of the first and second kinds, degenerate Bernoulli numbers of the second kind and degenerate Bell numbers and polynomials. Especially, it turns out that the degenerate Stirling numbers of the first and second kind appear very frequently when we study degenerate versions of some special numbers and polynomials. It is noteworthy that studying degenerate versions is not only limited to polynomials but also extended to transcendental functions. Indeed, the degenerate gamma functions were introduced in connection with degenerate Laplace transforms in [12]. It is also remarkable that the degenerate umbral calculus is introduced as a degenerate version of the classical umbral calculus. Indeed, the Sheffer sequences occupy the central position in the umbral calculus and are characterized by the generating functions for Sheffer pairs where the usual exponential function enters (see [8]). The motivation for the paper [8] started from the question that what if the usual exponential function is replaced by the degenerate exponential functions in \eqref{3}. With this replacement the $\lambda$-Sheffer polynomials are introduced, which are charactered by the generating functions for Sheffer pairs where the degenerate exponential function enters. \par
The Stirling number of the second $S_{2}(n,k)$ is the number of ways to partition a set of $n$ objects into $k$ nonempty subsets (see \eqref{2}). The (signed) Stirling number of the first kind $S_{1}(n,k)$ is defined such that the number of permutations of $n$ elements having exactly $k$ cycles is the nonnegative integer $(-1)^{n-k}S_{1}(n,k)=|S_{1}(n,k)|$ (see \eqref{1}). The degenerate Stirling numbers of the second kind $S_{2,\lambda}(n,k)$ (see \eqref{11}, \eqref{12}) and of the first kind $S_{1,\lambda}(n,k)$ (see \eqref{10}, \eqref{13}) appear most naturally when we replace the powers of $x$ by the generalized falling factorial polynomials $(x)_{k,\lambda}$ in the defining equations (see \eqref{1}, \eqref{2}, \eqref{10}, \eqref{11}).\par

Degenerate versions of special numbers and polynomials have been explored by various methods, including combinatorial methods, generating functions, umbral calculus techniques, $p$-adic analysis, differential equations, special functions, probability theory and analytic number theory. The aim of this paper is by using generating functions to further study some identities and properties on the degenerate Stirling numbers of the second kind, the degenerate $r$-Stirling numbers of the second kind, the degenerate Stirling numbers of the first kind and the unsigned degenerate Stirling numbers of the first kind. \par

The outline of this paper is as follows. In Section 1, we recall the facts that are needed throughout this paper.
The main results are obtained in Section 2.
Some identities on the degenerate Stirling numbers of the second kind are derived in Theorems 1, 10, 13 and 14,
those on the $r$-Stirling numbers of the second kind are stated in Theorem 2, those on the degenerate Stirling numbers of the first kind are obtained in Theorem 4, 5, 7 and 8. Finally, certain identities on the unsigned degenerate Stirling numbers of the first kind are obtained in Theorems 3, 6, 9 and 15.

For $n\ge 0$, the Stirling numbers of the first kind are defined by
\begin{equation}
(x)_{n}=\sum_{k=0}^{n}S_{1}(n,k)x^{k},\quad (\mathrm{see}\ [1-20]),\label{1}	
\end{equation}
where $(x)_{0}=1$, $(x)_{n}=x(x-1)\cdots(x-(n-1)),\ (n\ge 1)$. \par
As the inversion formula of \eqref{1} the Stirling numbers of the second kind are defined by
\begin{equation}
x^{n}=\sum_{k=0}^{n}S_{2}(n,k)(x)_{k},\quad(\mathrm{see}\ [2,5,14,17]).\label{2}	
\end{equation}
For any $\lambda\in\mathbb{R}$, the degenerate exponential functions are defined by
\begin{equation}
e_{\lambda}^{x}(t)=\sum_{k=0}^{\infty}\frac{(x)_{k,\lambda}}{k!}t^{k},\quad e_{\lambda}(t)=e_{\lambda}^{1}(t),\quad (\mathrm{see}\ [12,14]),\label{3}
\end{equation}
where the generalized falling factorials $(x)_{n,\lambda}$ are given by \\
$(x)_{0,\lambda}=1,\ (x)_{n,\lambda}=x(x-\lambda)\cdots(x-(n-1)\lambda),\ (n\ge 1)$. Note that $\displaystyle\lim_{\lambda\rightarrow 0}e_{\lambda}^{x}(t)=e^{xt}\displaystyle$. \par
For $k\ge 0$, the Lah numbers are given by
\begin{equation}
\frac{1}{k!}\bigg(\frac{t}{1-t}\bigg)^{k}=\sum_{n=k}^{\infty}L(n,k)\frac{t^{n}}{n!},\quad (\mathrm{see}\ [2,5,14,17]).\label{4}	
\end{equation}
From \eqref{4}, we note that
\begin{equation}
\langle x\rangle_{n}=\sum_{k=0}^{n}L(n,k)(x)_{k},\quad (n\ge 0),\quad (\mathrm{see}\ [2,5,14,17]),\label{5}
\end{equation}
where $\langle x\rangle_{0}=1,\ \langle x\rangle_{n}=x(x+1)\cdots(x+n-1),\ (n\ge 1)$. \par
By \eqref{5}, we get
\begin{equation}
L(n,k)=\frac{n!}{k!}\binom{n-1}{k-1},\quad (n,k\ge 1),\quad (\mathrm{see}\ [2,7,10,11]).\label{6}
\end{equation}
For any $\alpha\in\mathbb{R}$, the generalized Laguerre polynomials $L_{n}^{(\alpha)}(x)$ are given by
\begin{equation}
\begin{aligned}
	L_{n}^{(\alpha)}(x)&=\frac{1}{n!}x^{-\alpha}e^{x}\frac{d^{n}}{dx^{n}}\Big(e^{-x}x^{n+\alpha}\Big) \\
	&=x^{-\alpha}\frac{1}{n!}\bigg(\frac{d}{dx}-1\bigg)^{n}x^{n+\alpha},\quad (\mathrm{see}\ [2,7,17]).
\end{aligned}\label{7}
\end{equation}
When $\alpha=0$, $L_{n}(x)=L_{n}^{(0)}(x)$ are called the Laguerre polynomials. \par
Let $\log_{\lambda}t$ be the compositional inverse function of $e_{\lambda}(t)$. Then we have
\begin{equation}
\log_{\lambda}(1+t)=\sum_{n=1}^{\infty}\lambda^{n-1}(1)_{n,1/\lambda}\frac{t^{n}}{n!},\quad (\mathrm{see}\ [7]). \label{9}	
\end{equation}
We note that $\displaystyle\lim_{\lambda\rightarrow 0}\log_{\lambda}(1+t)=\log(1+t)\displaystyle$.\par
It is known that the degenerate Stirling numbers of the first kind are defined by
\begin{equation}
(x)_{n}=\sum_{k=0}^{n}S_{1,\lambda}(n,k)(x)_{k,\lambda},\quad (n\ge 0),\quad (\mathrm{see}\ [7]).\label{10}
\end{equation}
The degenerate Stirling numbers of the second kind are given by
\begin{equation}
(x)_{n,\lambda}=\sum_{k=0}^{n}S_{2,\lambda}(n,k)(x)_{k},\quad (n\ge 0),\quad (\mathrm{see}\ [7]). \label{11}
\end{equation}
From \eqref{10} and \eqref{11}, we note that $\displaystyle\lim_{\lambda\rightarrow 0}S_{1,\lambda}(n,k)=S_{1}(n,k)\displaystyle$ and  $\displaystyle\lim_{\lambda\rightarrow 0}S_{2,\lambda}(n,k)=S_{2}(n,k),\ (n,k\ge 0)\displaystyle$.\par
From [13], we recall that
\begin{equation}
\frac{1}{k!}\Big(e_{\lambda}(t)-1\Big)^{k}=\sum_{n=k}^{\infty}S_{2,\lambda}(n,k)\frac{t^{n}}{n!}, \label{12}
\end{equation}
and
\begin{equation}
\frac{1}{k!}\Big(\log_{\lambda}(1+t)\Big)^{k}=\sum_{n=k}^{\infty}S_{1,\lambda}(n,k)\frac{t^{n}}{n!}\quad (k\ge 0).\label{13}	
\end{equation}

\section{Some identities and properties on degenerate Stirling numbers}
For $k\in\mathbb{Z}$ with $k\ge 0$, we have
\begin{equation}
\begin{aligned}
e_{\lambda}(t)\frac{1}{k!}\big(e_{\lambda}(t)-1\big)^{k}&=\sum_{l=0}^{\infty}\frac{(1)_{l,\lambda}}{l!}t^{l}\sum_{j=k}^{\infty}S_{2,\lambda}(j,k)\frac{t^{j}}{j!}\\
&=\sum_{n=k}^{\infty}\bigg(\sum_{j=k}^{n}\binom{n}{j}S_{2,\lambda}(j,k)(1)_{n-j,\lambda}\bigg)\frac{t^{n}}{n!}.
\end{aligned}	\label{14}
\end{equation}
On the other hand,
\begin{equation}
\begin{aligned}
&e_{\lambda}(t)\frac{1}{k!}\big(e_{\lambda}(t)-1\big)^{k}=\big(e_{\lambda}(t)-1+1\big)\frac{1}{k!}\big(e_{\lambda}(t)-1\big)^{k} \\
&=\frac{k+1}{(k+1)!}\big(e_{\lambda}(t)-1\big)^{k+1}+\frac{1}{k!}\big(e_{\lambda}(t)-1\big)^{k}=\sum_{n=k}^{\infty}\Big((k+1)S_{2,\lambda}(n,k+1)+S_{2,\lambda}(n,k)\Big)\frac{t^{n}}{n!}.
\end{aligned}\label{15}
\end{equation}
Thus, by \eqref{14} and \eqref{15}, we get
\begin{equation}
\sum_{j=k}^{n}\binom{n}{j}S_{2,\lambda}(j,k)(1)_{n-j,\lambda}=(k+1)S_{2,\lambda}(n,k+1)+S_{2,\lambda}(n,k).\label{16}	
\end{equation}
From \eqref{11}, we note that
\begin{align}
\sum_{k=0}^{n+1}S_{2,\lambda}(n+1,k)(x)_{k}&=(x)_{n+1,\lambda}=(x)_{n,\lambda}(x-n\lambda)\label{17} \\
&=x\sum_{k=0}^{n}S_{2,\lambda}(n,k)(x)_{k}-n\lambda \sum_{k=0}^{n}S_{2,\lambda}(n,k)(x)_{k}\nonumber\\
&=\sum_{k=0}^{n}S_{2,\lambda}(n,k)(x)_{k}(x-k+k)-n\lambda\sum_{k=0}^{n}S_{2,\lambda}(n,k)(x)_{k}\nonumber \\
&=\sum_{k=0}^{n}S_{2,\lambda}(n,k)(x)_{k+1}+\sum_{k=0}^{n}(k-n\lambda)S_{2,\lambda}(n,k)(x)_{k}\nonumber \\
&=\sum_{k=1}^{n+1}S_{2,\lambda}(n,k-1)(x)_{k}+\sum_{n=0}^{n}(k-n\lambda)S_{2,\lambda}(n,k)(x)_{k}\nonumber \\
&=\sum_{k=0}^{n+1}\Big(S_{2,\lambda}(n,k-1)+(k-n\lambda)S_{2,\lambda}(n,k)\Big)(x)_{k}.\nonumber
\end{align}
Comparing the coefficients on both sides of \eqref{7}, we have
\begin{equation}
S_{2,\lambda}(n+1,k)=S_{2,\lambda}(n,k-1)+(k-n\lambda)S_{2,\lambda}(n,k).\label{18}
\end{equation}
Replacing $k$ by $k+1$, we get
\begin{equation}
S_{2,\lambda}(n+1,k+1)+n\lambda	 S_{2,\lambda}(n,k+1)=S_{2,\lambda}(n,k)+(k+1)S_{2,\lambda}(n,k+1). \label{19}
\end{equation}
Therefore, by \eqref{16} and \eqref{19}, we obtain the following theorem.
\begin{theorem}
For $n,k\in\mathbb{Z}$ with $n \ge k \ge 0$, we have
\begin{displaymath}
\sum_{j=k}^{n}\binom{n}{j}S_{2,\lambda}(j,k)(1)_{n-j,\lambda}=S_{2,\lambda}(n+1,k+1)+n\lambda S_{2,\lambda}(n,k+1).
\end{displaymath}	
\end{theorem}
It is known that the degenerate $r$-Stirling numbers of the second kind are defined by
\begin{equation}
e_{\lambda}^{r}(t)\frac{1}{k!}\big(e_{\lambda}(t)-1\big)^{k}=\sum_{n=k}^{\infty}S_{2,\lambda}^{(r)}(n+r,k+r)\frac{t^{n}}{n!},\quad (\mathrm{see}\ [10,11]). \label{20}
\end{equation}
From \eqref{20}, we note that
\begin{equation}
\begin{aligned}
\frac{1}{k!}\big(e_{\lambda}(t)-1\big)^{k}e_{\lambda}^{r}(t)&=\sum_{l=k}^{\infty}S_{2,\lambda}(l,k)\frac{t^{l}}{l!}\sum_{m=0}^{\infty}(r)_{m,\lambda}\frac{t^{m}}{m!} \\
&=\sum_{n=k}^{\infty}\bigg(\sum_{l=k}^{n}\binom{n}{l}S_{2,\lambda}(l,k)(r)_{n-l,\lambda}\bigg)\frac{t^{n}}{n!}.\end{aligned}\label{21}	
\end{equation}
By \eqref{20} and \eqref{21}, we get
\begin{equation}
S_{2,\lambda}^{(r)}(n+r,k+r)=\sum_{l=k}^{n}\binom{n}{l}S_{2,\lambda}(l,k)(r)_{n-l,\lambda}\label{22},
\end{equation}
where $n,k,r$ are nonnegative integers. \par
From \eqref{21}, we note that
\begin{align}
\sum_{n=0}^{\infty}\bigg(\sum_{k=0}^{n}S_{2,\lambda}^{(r)}(n+r,k+r)(x)_{k}\bigg)\frac{t^{n}}{n!}&=\sum_{k=0}^{\infty}\bigg(\sum_{n=k}^{\infty}S_{2,\lambda}^{(r)}(n+r,k+r)\frac{t^{n}}{n!}\bigg)(x)_{k} \label{23} \\
&=\sum_{k=0}^{\infty}\bigg(\frac{1}{k!}\Big(e_{\lambda}(t)-1\Big)^{k}e_{\lambda}^{r}(t)\bigg)(x)_{k}\nonumber \\
&=e_{\lambda}^{r}(t)\sum_{k=0}^{\infty}(x)_{k}\frac{1}{k!}\Big(e_{\lambda}(t)-1\Big)^{k}\nonumber \\
&=e_{\lambda}^{r}(t)e_{\lambda}^{x}(t)=\sum_{n=0}^{\infty}(x+r)_{n,\lambda}\frac{t^{n}}{n!}. \nonumber
\end{align}
Comparing the coefficients on both sides of \eqref{23}, we have
\begin{align}
(x+r)_{n,\lambda}=\sum_{k=0}^{n}S_{2,\lambda}^{(r)}(n+r,k+r)(x)_{k},\quad (\mathrm{see}\ [14]).\label{23-1}
\end{align}
Further, we note that
\begin{align}
(x+r)_{n,\lambda}&=\sum_{l=0}^{n}S_{2,\lambda}(n,l)(x+r)_{l}\nonumber\\
&=\sum_{l=0}^{n}S_{2,\lambda}(n,l)\sum_{k=0}^{l}\binom{l}{k}(x)_{k}(r)_{l-k}\label{23-2}\\
&=\sum_{k=0}^{n}\sum_{l=k}^{n}\binom{l}{k}S_{2,\lambda}(n,l)(r)_{l-k}(x)_{k}.\nonumber
\end{align}
From \eqref{22}, \eqref{23}, \eqref{23-1} and \eqref{23-2}, we obtain the next theorem.
\begin{theorem}
For $n,k,r\in\mathbb{Z}$ with $n \ge k \ge 0$, and $r \ge 0$, we have
\begin{align*}
&S_{2,\lambda}^{(r)}(n+r,k+r)=\sum_{l=k}^{n}\binom{n}{l}S_{2,\lambda}(l,k)(r)_{n-l,\lambda}=\sum_{l=k}^{n}\binom{l}{k}S_{2,\lambda}(n,l)(r)_{l-k},\\
&(x+r)_{n,\lambda}=\sum_{k=0}^{n}S_{2,\lambda}^{(r)}(n+r,k+r)(x)_{k}.\\
\end{align*}
\end{theorem}
It is known that the unsigned degenerate Stirling numbers of the first kind are defined by
\begin{equation}
\langle x\rangle_{n}=\sum_{k=0}^{n}{n \brack k}_{\lambda}\langle x\rangle_{k,\lambda},\quad (n\ge 0),\quad (\mathrm{see}\ [14]).\label{24}	
\end{equation}
By \eqref{10} and \eqref{24}, we easily get
\begin{equation}
{n \brack k}_{\lambda}=(-1)^{n-k}S_{1,\lambda}(n,k),\quad (n,k\ge 0). \label{25}	
\end{equation}
Now, we observe that
\begin{align}
\frac{(-1)^{k}}{k!}\big(\log_{\lambda}(1-t)\big)^{k}\frac{1}{1-t}&=\frac{(-1)^{k}}{k!}\big(\log_{\lambda}(1-t)\big)^{k}e_{\lambda}^{-1}\big(\log_{\lambda}(1-t)\big) \label{26}\\
&=\frac{(-1)^{k}}{k!}\big(\log_{\lambda}(1-t)\big)^{k}\sum_{l=0}^{\infty}(-1)^{l}\langle 1\rangle_{l,\lambda}\frac{1}{l!}\big(\log_{\lambda}(1-t)\big)^{l}\nonumber \\
&=\frac{1}{k!}\sum_{l=0}^{\infty}(-1)^{k+l}\langle 1\rangle_{l,\lambda}\frac{1}{l!}\big(\log_{\lambda}(1-t)\big)^{l+k}\nonumber 	\\
&=\frac{1}{k!}\sum_{l=k}^{\infty}(-1)^{l}\langle 1\rangle_{l-k,\lambda}\frac{l!}{(l-k)!}\frac{1}{l!}\big(\log_{\lambda}(1-t)\big)^{l}\nonumber\\
&=\sum_{l=k}^{\infty}(-1)^{l}\langle 1\rangle_{l-k,\lambda}\binom{l}{k}\sum_{n=l}^{\infty}S_{1,\lambda}(n,l)(-1)^{n}\frac{t^{n}}{n!} \nonumber \\
&=\sum_{n=k}^{\infty}\bigg(\sum_{l=k}^{n}(-1)^{n-l}\binom{l}{k}S_{1,\lambda}(n,l)\langle 1\rangle_{l-k,\lambda}\bigg)\frac{t^{n}}{n!}. \nonumber
\end{align}
On the other hand,
\begin{align}
\frac{(-1)^{k}}{k!}\big(\log_{\lambda}(1-t)\big)^{k}\frac{1}{1-t}&=(-1)^{k}\sum_{l=k}^{\infty}S_{1,\lambda}(l,k)(-1)^{l}\frac{t^{l}}{l!}\sum_{m=0}^{\infty}t^{m} \label{27} \\
&=\sum_{n=k}^{\infty}\Big(\sum_{l=k}^{n}\frac{S_{1\lambda}(l,k)}{l!}(-1)^{l-k}\Big)t^{n}\nonumber \\
&=\sum_{n=k}^{\infty}\Big(n!\sum_{l=k}^{n}(-1)^{l-k}\frac{S_{1,\lambda}(l,k)}{l!}\Big)\frac{t^{n}}{n!}.\nonumber	
\end{align}
Therefore, by \eqref{26} and \eqref{27}, we obtain the following theorem.
\begin{theorem}
For $n,k\in\mathbb{Z}$ with $n \ge k\ge 0$, we have
\begin{displaymath}
\frac{1}{n!}\sum_{l=k}^{n}\binom{l}{k}\langle 1\rangle_{l-k,\lambda}{n \brack l}_{\lambda}=\sum_{l=k}^{n}\frac{1}{l!}{l \brack k}_{\lambda}.
\end{displaymath}	
\end{theorem}
From \eqref{13}, we note that
\begin{align}
\frac{1}{k!}\big(\log_{\lambda}(1+t)\big)^{k}(1+t)&=\sum_{n=k}^{\infty}S_{1,\lambda}(n,k)\frac{t^{n}}{n!}(1+t)\label{28} \\
&=\sum_{n=k}^{\infty}S_{1,\lambda}(n,k)\frac{t^{n}}{n!}+\sum_{n=k+1}^{\infty}nS_{1,\lambda}(n-1,k)\frac{t^{n}}{n!}\nonumber\\
&=\sum_{n=k}^{\infty}\big(S_{1,\lambda}(n,k)+nS_{1,\lambda}(n-1,k)\big)\frac{t^{n}}{n!}.\nonumber
\end{align}
On the other hand
\begin{align}
\frac{1}{k!}\big(\log_{\lambda}(1+t)\big)^{k}e_{\lambda}\big(\log_{\lambda}(1+t)\big)&=\frac{1}{k!}\sum_{l=0}^{\infty}\frac{(1)_{l,\lambda}}{l!}\big(\log_{\lambda}(1+t)\big)^{k+l} \label{29} 	\\
&=\frac{1}{k!}\sum_{l=k}^{\infty}\frac{(1)_{l-k,\lambda}l!}{(l-k)!}\frac{1}{l!}\big(\log_{\lambda}(1+t)\big)^{l}\nonumber \\
&=\sum_{n=k}^{\infty}\bigg(\sum_{l=k}^{n}\binom{l}{k}(1)_{l-k,\lambda}S_{1,\lambda}(n,l)\bigg)\frac{t^{n}}{n!}.\nonumber
\end{align}
Therefore, by \eqref{28} and \eqref{29}, we obtain the following theorem.
\begin{theorem}
For $n,k\in\mathbb{N}$ with $n \ge k \ge 0$, we have
\begin{displaymath}
S_{1,\lambda}(n,k)+nS_{1,\lambda}(n-1,k)=\sum_{l=k}^{n}\binom{l}{k}S_{1,\lambda}(n,l)(1)_{l-k,\lambda}.
\end{displaymath}	
\end{theorem}
For any $\alpha\in\mathbb{R}$, and $k\in\mathbb{Z}$ with $k\ge 0$, we have
\begin{align}
\frac{1}{k!}(1+t)^{\alpha}\big(\log_{\lambda}(1+t)\big)^{k}&=\frac{1}{k!}\sum_{l=0}^{\infty}(\alpha)_{l,\lambda}\frac{1}{l!}\big(\log_{\lambda}(1+t)\big)^{l+k}\label{30} \\
&=\sum_{l=k}^{\infty}\binom{l}{k}(\alpha)_{l-k,\lambda}\sum_{n=l}^{\infty}S_{1,\lambda}(n,l)\frac{t^{n}}{n!}\nonumber \\
&=\sum_{n=k}^{\infty}\bigg(\sum_{l=k}^{n}\binom{l}{k}(\alpha)_{n-k,\lambda}S_{1,\lambda}(n,l)\bigg)\frac{t^{n}}{n!}. 	\nonumber
\end{align}
By binomial expansion, we get
\begin{align}
\frac{1}{k!}(1+t)^{\alpha}\big(\log_{\lambda}(1+t)\big)^{k}&=\sum_{m=0}^{\infty}\binom{\alpha}{m}t^{m}\sum_{l=k}^{\infty}S_{1,\lambda}(l,k)\frac{t^{l}}{l!} \label{31}	\\
&=\sum_{n=k}^{\infty}\sum_{l=k}^{n}\frac{S_{1,\lambda}(l,k)}{l!}\binom{\alpha}{n-l}n!\frac{t^{n}}{n!}\nonumber\\
&=\sum_{n=k}^{\infty}\bigg(\sum_{l=k}^{n}\binom{n}{l}\binom{\alpha}{n-l}S_{1,\lambda}(l,k)(n-l)!\bigg)\frac{t^{n}}{n!}.\nonumber
\end{align}
Therefore, by \eqref{30} and \eqref{31}, we obtain the following theorem.
\begin{theorem}
For $\alpha\in\mathbb{R}$ and $n,k\in\mathbb{Z}$ with $n \ge k \ge 0$, we have
\begin{align*}
\sum_{l=k}^{n}\binom{l}{k}(\alpha)_{l-k,\lambda}S_{1,\lambda}(n,l)&=\sum_{l=k}^{n}\binom{n}{l}\binom{\alpha}{n-l}S_{1,\lambda}(l,k)(n-l)!\\
&=\sum_{l=0}^{n-k}\binom{n}{l}\binom{\alpha}{l}l!S_{1,\lambda}(n-l,k). 	
\end{align*}
\end{theorem}
From \eqref{13}, we have
\begin{align}
\frac{1}{t}\big(\log_{\lambda}(1+t)\big)^{k}&=\frac{k!}{t}\frac{1}{k!}\big(\log_{\lambda}(1+t)\big)^{k}\label{32}\\
&=\frac{k!}{t}\sum_{n=k}^{\infty}S_{1,\lambda}(n,k)\frac{t^{n}}{n!} \nonumber \\
&=k!\sum_{n=k-1}^{\infty}S_{1,\lambda}(n+1,k)\frac{t^{n}}{(n+1)!}.\nonumber	
\end{align}
Thus, by \eqref{32}, we get
\begin{align}
\frac{1}{1-t}\frac{1-t}{t}\Big(\log_{\lambda}\Big(1+\frac{t}{1-t}\Big)\Big)^{k}&=\frac{1}{1-t}k!\sum_{l=k-1}^{\infty}S_{1,\lambda}(l+1,k)\frac{\big(\frac{t}{1-t}\big)^{l}}{(l+1)!} \label{33} \\
&=k!\sum_{l=k-1}^{\infty}S_{1,\lambda}(l+1,k)\frac{t^{l}}{(l+1)!}\Big(\frac{1}{1-t}\Big)^{l+1}\nonumber\\
&=k!\sum_{l=k-1}^{\infty}S_{1,\lambda}(l+1,k)\frac{t^{l}}{(l+1)!}\sum_{m=0}^{\infty}\binom{l+m}{m}t^{m}\nonumber\\
&=k!\sum_{n=k-1}^{\infty}\bigg(\sum_{l=k-1}^{n}\frac{S_{1,\lambda}(l+1,k)}{(l+1)!}\binom{n}{l}\bigg)t^{n}.\nonumber	
\end{align}
On the other hand, by \eqref{32}, we see that
\begin{align}
&\frac{1}{1-t}\frac{1-t}{t}\Big(\log_{\lambda}\Big(1+\frac{t}{1-t}\Big)\Big)^{k}=\frac{1}{t}\Big(\log_{\lambda}\frac{1}{1-t}\Big)^{k}\label{34}\\
&=\frac{1}{t}\big(-\log_{-\lambda}(1-t)\big)^{k}=(-1)^{k-1}\frac{1}{-t}\big(\log_{-\lambda}(1-t)\big)^{k}\nonumber	\\
&=(-1)^{k-1}k!\sum_{n=k-1}^{\infty}(-1)^{n}\frac{S_{1,-\lambda}(n+1,k)}{(n+1)!}t^{n}.\nonumber
\end{align}
Therefore, by \eqref{33} and \eqref{34}, we obtain the following theorem.
\begin{theorem}
For $n,k\in\mathbb{N}$ with $n\ge k-1$,  we have
\begin{displaymath}
\frac{1}{(n+1)!}{n+1 \brack k}_{-\lambda}=\sum_{l=k-1}^{n}\binom{n}{l}\frac{S_{1,\lambda}(l+1,k)}{(l+1)!}.
\end{displaymath}
\end{theorem}
We observe that
\begin{align}
\frac{1}{t}\big(\log_{\lambda}(1+t)\big)^{k+1}&=k!\log_{\lambda}(1+t)\frac{1}{t}\frac{1}{k!}\big(\log_{\lambda}(1+t)\big)^{k}\label{35}\\
&=k!\sum_{j=1}^{\infty}\lambda^{j-1}(1)_{j,1/\lambda}\frac{t^{j}}{j!}\sum_{l=k-1}^{\infty}\frac{1}{(l+1)!}S_{1,\lambda}(l+1,k)t^{l}\nonumber \\
&=\sum_{n=k}^{\infty}\bigg(k!\sum_{l=k-1}^{n-1}\lambda^{n-l-1}(1)_{n-l,1/\lambda}\frac{1}{(n-l)!}\frac{1}{(l+1)!}S_{1,\lambda}(l+1,k)\bigg)t^{n}.\nonumber
\end{align}
From \eqref{13}, we also note that
\begin{align}
\frac{1}{t}\big(\log_{\lambda}(1+t)\big)^{k+1}&=\frac{(k+1)!}{t}\frac{1}{(k+1)!}\big(\log_{\lambda}(1+t)\big)^{k+1}\label{36} \\
&=\frac{(k+1)!}{t}\sum_{n=k+1}^{\infty}S_{1,\lambda}(n,k+1)\frac{t^{n}}{n!}\nonumber\\
&=(k+1)!\sum_{n=k}^{\infty}\frac{S_{1,\lambda}(n+1,k+1)}{(n+1)!}t^{n}.\nonumber
\end{align}
Therefore, by \eqref{35} and \eqref{36}, we obtain the following theorem.
\begin{theorem}
For $n,k\in\mathbb{Z}$ with $n \ge k\ge 0$, we have
\begin{displaymath}
S_{1,\lambda}(n+1,k+1)=\frac{\lambda^{n-1}}{k+1}(n+1)!\sum_{l=k-1}^{n-1}\frac{\lambda^{-l}}{(n-l)!}(1)_{n-l,1/\lambda}\frac{S_{1,\lambda}(l+1,k)}{(l+1)!}.
\end{displaymath}
\end{theorem}
For $n,k\in\mathbb{N}$, by \eqref{10}, we easily get
\begin{equation}
S_{1,\lambda}(n+1,k)=S_{1,\lambda}(n,k-1)+(k\lambda-n)S_{1,\lambda}(n,k).\label{37}
\end{equation}
From \eqref{13}, we note that
\begin{align}
\big(\log_{\lambda}(1+t)\big)^{p+1}&=(p+1)!\frac{1}{(p+1)!}\big(\log_{\lambda}(1+t)\big)^{p+1}\label{38} \\
&=(p+1)!\sum_{n=p+1}^{\infty}S_{1,\lambda}(n,p+1)\frac{t^{n}}{n!}. \nonumber
\end{align}
Taking derivatives on both sides of \eqref{38}, we have
\begin{equation}
(p+1)\frac{1}{1+t}\big(\log_{\lambda}(1+t)\big)^{p}(1+t)^{\lambda}=(p+1)!\sum_{n=p}^{\infty}S_{1,\lambda}(n+1,p+1)\frac{t^{n}}{n!}. \label{39}
\end{equation}
Thus, by \eqref{30}, we get
\begin{align}
\sum_{n=p}^{\infty}S_{1,\lambda}(n+1,p+1)\frac{t^{n}}{n!}&=\frac{1}{p!}\frac{1}{1+t}\big(\log_{\lambda}(1+t)\big)^{p}(1+t)^{\lambda}\label{40} \\
&=\sum_{j=0}^{\infty}(-1)^{j}t^{j}\sum_{l=p}^{\infty}S_{1,\lambda}(l,p)\frac{t^{l}}{l!}\sum_{m=0}^{\infty}(\lambda)_{m}\frac{t^{m}}{m!}\nonumber \\
&=\sum_{k=p}^{\infty}\sum_{l=p}^{k}(-1)^{k-l}\frac{S_{1,\lambda}(l,p)}{l!}t^{k}\sum_{m=0}^{\infty}(\lambda)_{m}\frac{t^{m}}{m!}\nonumber \\
&=\sum_{n=p}^{\infty}\bigg(\sum_{k=p}^{n}\sum_{l=p}^{k}(-1)^{k-l}\binom{n}{k}\frac{S_{1,\lambda}(l,p)}{l!}(\lambda)_{n-k}k!\bigg)\frac{t^{n}}{n!}. \nonumber
\end{align}
Comparing the coefficients on both sides of \eqref{40}, we obtain the following theorem.
\begin{theorem}
For $n,p\in\mathbb{Z}$ with $n \ge p\ge 0$, we have
\begin{displaymath}
S_{1,\lambda}(n+1,p+1)=\sum_{k=p}^{n}\sum_{l=p}^{k}(-1)^{k-l}\binom{k}{l}\binom{n}{k}(k-l)!S_{1,\lambda}(l,p)(\lambda)_{n-k}.
\end{displaymath}	
\end{theorem}
Note that
\begin{displaymath}
S_{1}(n+1,p+1)=\lim_{\lambda\rightarrow 0}S_{1,\lambda}(n+1,p+1)=\sum_{l=p}^{n}(-1)^{n-l}\binom{n}{l}(n-l)!S_{1}(l,p).
\end{displaymath}
For $p\in\mathbb{Z}$ with $p\ge 0$, we have
\begin{align}
&\frac{1}{1-t}\Big(\log_{\lambda}\Big(1+\frac{t}{1-t}\Big)\Big)^{p}=(-1)^{p}\frac{p!}{1-t}\frac{1}{p!}\Big(\log_{-\lambda}(1-t)\Big)^{p}\label{41} \\
&=(-1)^{p}\frac{p!}{1-t}\sum_{k=p}^{\infty}S_{1,-\lambda}(k,p)(-1)^{k}\frac{t^{k}}{k!}\nonumber \\
&=(-1)^{p}p!\sum_{m=0}^{\infty}t^{m}\sum_{k=p}^{\infty}S_{1,-\lambda}(k,p)(-1)^{k}\frac{t^{k}}{k!}\nonumber \\
&=\sum_{n=p}^{\infty}\Big((-1)^{p}p!\sum_{k=p}^{n}\frac{S_{1,-\lambda}(k,p)}{k!}(-1)^{k}\Big)t^{n}.\nonumber	
\end{align}
On the other hand, by \eqref{13}, we get
\begin{align}
&\frac{1}{1-t}\Big(\log_{\lambda}\Big(1+\frac{t}{1-t}\Big)\Big)^{p}=p!\sum_{k=p}^{\infty}S_{1,\lambda}(k,p)\frac{t^{k}}{k!}\Big(\frac{1}{1-t}\Big)^{k+1}\label{42} \\
&=p!\sum_{k=p}^{\infty}S_{1,\lambda}(k,p)\frac{t^{k}}{k!}\sum_{m=0}^{\infty}\binom{k+m}{m}t^{m} \nonumber \\
&=p!\sum_{n=p}^{\infty}\bigg(\sum_{k=p}^{n}S_{1,\lambda}(k,p)\frac{1}{k!}\binom{n}{k}\bigg)t^{n}.\nonumber
\end{align}
Therefore, by \eqref{41} and \eqref{42}, we obtain the following theorem.
\begin{theorem}
For $n,p\in\mathbb{Z}$ with $n \ge p\ge 0$, we have
\begin{displaymath}
\sum_{k=p}^{n}\frac{1}{k!}{k \brack p}_{-\lambda}=\sum_{k=p}^{n}\frac{1}{k!}\binom{n}{k}S_{1,\lambda}(k,p).
\end{displaymath}	
\end{theorem}
It is known that the degenerate Bell polynomials are defined by
\begin{equation}
e^{x(e_{\lambda}(t)-1)}=\sum_{n=0}^{\infty}\phi_{n,\lambda}(x)\frac{t^{n}}{n!},\quad (\mathrm{see}\ [13]).\label{43}
\end{equation}
From \eqref{43}, we note that
\begin{equation}
\phi_{n,\lambda}(x)=\sum_{k=0}^{n}S_{2,\lambda}(n,k)x^{k},\quad (n\ge 0). \label{44}	
\end{equation}
By \eqref{3} and \eqref{13}, we get
\begin{align}
&e_{\lambda}(t)\big(\log_{\lambda}(1+t)\big)^{p}=p!e_{\lambda}(t)\frac{1}{p!}\big(\log_{\lambda}(1+t)\big)^{p} \nonumber \\
&=p!\sum_{l=0}^{\infty}(1)_{l,\lambda}\frac{t^{l}}{l!}\sum_{k=p}^{\infty}S_{1,\lambda}(k,p)\frac{t^{k}}{k!}=p!\sum_{n=p}^{\infty}\bigg(\sum_{k=p}^{n}\binom{n}{k}S_{1,\lambda}(k,p)(1)_{n-k,\lambda}\bigg)\frac{t^{n}}{n!}.\nonumber	
\end{align}
From \eqref{12}, we note that
\begin{align}
&x^{p}\frac{1}{p!}\big(e_{\lambda}(t)-1\big)^{p}e^{x(e_{\lambda}(t)-1)}=\frac{1}{p!}\sum_{k=0}^{\infty}x^{k+p}\frac{1}{k!}\big(e_{\lambda}(t)-1\big)^{k+p}\label{45} \\
&=\sum_{k=p}^{\infty}\binom{k}{p}x^{k}\frac{1}{k!}\big(e_{\lambda}(t)-1\big)^{k}=\sum_{k=p}^{\infty}\binom{k}{p}x^{k}\sum_{n=k}^{\infty}S_{2,\lambda}(n,k)\frac{t^{n}}{n!}\nonumber \\
&=\sum_{n=p}^{\infty}\bigg(\sum_{k=p}^{n}\binom{k}{p}x^{k}S_{2,\lambda}(n,k)\bigg)\frac{t^{n}}{n!}.\nonumber
\end{align}
On the other hand,
\begin{align}
x^{p}\frac{1}{p!}(e_{\lambda}(t)-1)^{p}e^{x(e_{\lambda}(t)-1)}&=x^{p}\sum_{k=p}^{\infty}S_{2,\lambda}(k,p)\frac{t^{k}}{k!}\sum_{m=0}^{\infty}\phi_{m,\lambda}(x)\frac{t^{m}}{m!}\label{46}\\
&=x^{p}\sum_{n=p}^{\infty}\bigg(\sum_{k=p}^{n}\binom{n}{k}S_{2,\lambda}(k,p)\phi_{n-k,\lambda}(x)\bigg)\frac{t^{n}}{n!}.\nonumber
\end{align}
Therefore, by \eqref{45} and \eqref{46}, we obtain the following theorem.
\begin{theorem}
For $n,p\in\mathbb{Z}$	with $n \ge p\ge 0$, we have
\begin{displaymath}
x^{p}\sum_{k=p}^{n}\binom{n}{k}S_{2,\lambda}(k,p)\phi_{n-k,\lambda}(x)=\sum_{k=p}^{n}\binom{k}{p}S_{2,\lambda}(n,k)x^{k},
\end{displaymath}
In particular, for $x=1$,
\begin{displaymath}
\sum_{k=p}^{n}\binom{n}{k}S_{2,\lambda}(k,p)\phi_{n-k,\lambda}=\sum_{k=p}^{n}\binom{k}{p}S_{2,\lambda}(n,k),
\end{displaymath}
where $\phi_{n,\lambda}=\phi_{n,\lambda}(1)$ are the degenerate Bell numbers.
\end{theorem}
For $p=1$, we have
\begin{displaymath}
\sum_{k=1}^{n}kS_{2,\lambda}(n,k)=\sum_{k=0}^{n-1}\binom{n}{k}\phi_{k,\lambda}(1)_{n-k,\lambda},\quad (n\ge 1).
\end{displaymath}
In view of \eqref{7}, we may consider the generalized degenerate Laguerre polynomials    given by
\begin{equation}
\frac{1}{n!}x^{-\alpha}e_{\lambda}(x)\bigg(\frac{d}{dx}\bigg)^{n}\Big(x^{n+\alpha}e_{\lambda}^{-1}(x)\Big)=L_{n,\lambda}^{(\alpha)}(x),\quad (n\ge 0). \label{47}
\end{equation}
From \eqref{47}, we can derive the following equation
\begin{equation}
L_{n,\lambda}^{(\alpha)}(x)=\sum_{k=0}^{n}\binom{n+\alpha}{n-k}(-1)^{k}\langle 1\rangle_{k,\lambda}\frac{1}{k!}\bigg(\frac{x}{1+\lambda x}\bigg)^{k}.\label{48}
\end{equation}
Now, we observe that
\begin{align}
&\frac{1}{p!}\bigg(\frac{t}{1-t}\frac{x}{1+\lambda x}\bigg)^{p}e_{\lambda}\bigg(\frac{t}{1-t}\frac{x}{1+\lambda x}\bigg)=\frac{1}{p!}\bigg(\frac{t}{1-t}\frac{x}{1+\lambda x}\bigg)^{p}\sum_{k=0}^{\infty}(1)_{k,\lambda}\frac{\big(\frac{x}{1+\lambda x}\big)^{k}}{(1-t)^{k}}\frac{t^{k}}{k!} \label{49}\\
&=\frac{1}{p!}\sum_{k=0}^{\infty}\frac{(1)_{k,\lambda}}{k!}\frac{\big(\frac{x}{1+\lambda x}\big)^{k+p}}{(1-t)^{k+p}}t^{k+p}=\frac{1}{p!}\sum_{k=p}^{\infty}\frac{(1)_{k-p,\lambda}}{(k-p)!}\bigg(\frac{t}{1-t}\bigg)^{k}\bigg(\frac{x}{1+\lambda x}\bigg)^{k}\nonumber\\
&=\frac{1}{p!}\sum_{k=p}^{\infty}\frac{k!(1)_{k-p,\lambda}}{(k-p)!}\bigg(\frac{x}{1+\lambda x}\bigg)^{k}\frac{1}{k!}\bigg(\frac{t}{1-t}\bigg)^{k} \nonumber \\
&=\sum_{k=p}^{\infty}\binom{k}{p}(1)_{k-p,\lambda}\bigg(\frac{x}{1+\lambda x}\bigg)^{k}\sum_{n=k}^{\infty}L(n,k)\frac{t^{n}}{n!} \nonumber \\
&=\sum_{n=p}^{\infty}\bigg(\sum_{k=p}^{n}\binom{k}{p}(1)_{k-p,\lambda}L(n,k)\bigg(\frac{x}{1+\lambda x}\bigg)^{k}\bigg)\frac{t^{n}}{n!}.\nonumber
\end{align}
Here we note that
\begin{equation}
\frac{1}{(1-t)^{\alpha+1}}e_{\lambda}^{-1}\bigg(\frac{t}{1-t}\frac{x}{1+\lambda x}\bigg)=\sum_{n=0}^{\infty}L_{n,\lambda}^{(\alpha)}(x)t^{n}.\label{50}
\end{equation}
Indeed, from \eqref{48} we deduce that
\begin{align*}
\sum_{n=0}^{\infty}L_{n,\lambda}^{(\alpha)}(x)t^{n}&=\sum_{k=0}^{\infty}(-1)^{k}\langle 1\rangle_{k,\lambda}\frac{1}{k!}\bigg(\frac{x}{1+\lambda x}\bigg)^{k}\sum_{n=k}^{\infty}\binom{n+\alpha}{n-k}t^{n}\\
&=\sum_{k=0}^{\infty}(-1)_{k,\lambda}\frac{1}{k!}\bigg(\frac{t}{1-t}\frac{x}{1+\lambda x}\bigg)^{k}(1-t)^{k}\sum_{n=k}^{\infty}\binom{n+\alpha}{n-k}t^{n-k}\\
&=\sum_{k=0}^{\infty}(-1)_{k,\lambda}\frac{1}{k!}\bigg(\frac{t}{1-t}\frac{x}{1+\lambda x}\bigg)^{k}
(1-t)^{k}(1-t)^{-\alpha-k-1}\\
&=(1-t)^{\alpha-1}e_{\lambda}^{-1}\bigg(\frac{t}{1-t}\frac{x}{1+\lambda x}\bigg).
\end{align*}
From \eqref{50}, we note that
\begin{align}
&\frac{1}{p!}\bigg(\frac{t}{1-t}\frac{x}{1+\lambda x}\bigg)^{p}e_{\lambda}\bigg(\frac{t}{1-t}\frac{x}{1+\lambda x}\bigg) \label{51}\\
&=\frac{t^{p}}{p!}\bigg(\frac{x}{1+\lambda x}\bigg)^{p}\frac{1}{(1-t)^{p}}e_{-\lambda}^{-1}\bigg(\frac{t}{1-t}\bigg(-\frac{x}{1+\lambda x}\bigg)\bigg)\nonumber \\
&=\frac{t^{p}}{p!}\bigg(\frac{x}{1+\lambda x}\bigg)^{p}\sum_{m=0}^{\infty}t^{m}L_{m,-\lambda}^{(p-1)}(-x)\nonumber \\
&=\frac{1}{p!}\bigg(\frac{x}{1+\lambda x}\bigg)^{p}\sum_{n=p}^{\infty}t^{n}L_{n-p,-\lambda}^{(p-1)}(-x)\nonumber\\
&=\frac{1}{p!}\bigg(\frac{x}{1+\lambda x}\bigg)^{p}\sum_{n=p}^{\infty}n!L_{n-p,-\lambda}^{(p-1)}(-x)\frac{t^{n}}{n!}.\nonumber	
\end{align}
Therefore, by \eqref{49} and \eqref{51}, we obtain the following theorem.
\begin{theorem}
For $n,p\in\mathbb{Z}$ with $n\ge p\ge 0$, we have \begin{displaymath}
L_{n-p,-\lambda}^{(p-1)}(-x)=\frac{p!}{n!}\sum_{k=p}^{n}\binom{k}{p}(1)_{k-p,\lambda}L(n,k)\bigg(\frac{x}{1+\lambda x}\bigg)^{k-p}.
\end{displaymath}	
\end{theorem}
We are going to use the following lemma.
\begin{lemma}
For $n\ge 0$, we have
	\begin{displaymath}
		a_{n,\lambda}=\sum_{k=0}^{n}S_{2,\lambda}(n,k)b_{k,\lambda}\ \Longleftrightarrow\ b_{n,\lambda}=\sum_{k=0}^{n}S_{1,\lambda}(n,k)a_{k,\lambda}
	\end{displaymath}
\end{lemma}
From \eqref{5} and \eqref{24}, we note that
\begin{align}
\sum_{k=0}^{n}L(n,k)(x)_{k}&=\langle x\rangle_{n}=\sum_{k=0}^{n}{n \brack k}_{\lambda}(x)_{k,\lambda}\label{52}\\
&=\sum_{k=0}^{n}{n \brack k}_{\lambda}\sum_{l=0}^{k}S_{2,\lambda}(k,l)(x)_{l}	\nonumber \\
&=\sum_{l=0}^{n}\bigg(\sum_{k=l}^{n}{n \brack k}_{\lambda}S_{2,\lambda}(k,l)\bigg)(x)_{l}\nonumber.
\end{align}
Comparing the coefficients on both sides of \eqref{52}, we get
\begin{equation}
L(n,l)=\sum_{k=l}^{n}{n\brack k}_{\lambda}S_{2,\lambda}(k,l),\quad (n\ge l\ge 0).\label{53}
\end{equation}
From \eqref{12}, we note that
\begin{align}
\frac{1}{p!}\bigg(\frac{1-e_{\lambda}(t)}{e_{\lambda}(t)}\bigg)^{p}&=\frac{1}{p!}\big(e_{\lambda}^{-1}(t)-1\big)^{p}=\frac{1}{p!}\big(e_{-\lambda}(-t)-1\big)^{p}\label{54} \\
&=\sum_{n=p}^{\infty}S_{2,-\lambda}(n,p)(-1)^{n}\frac{t^{n}}{n!}.\nonumber	
\end{align}
On the other hand, by \eqref{4}, we get
\begin{align}
&\frac{1}{p!}\bigg(\frac{1-e_{\lambda}(t)}{e_{\lambda}(t)}
\bigg)^{p}=\frac{1}{p!}\bigg(\frac{1-e_{\lambda}(t)}{e_{\lambda}(t)-1+1}\bigg)^{p}=\frac{1}{p!}\bigg(\frac{1-e_{\lambda}(t)}{1-(1-e_{\lambda}(t))}\bigg)^{p}\label{55} \\
&=\sum_{k=p}^{\infty}L(k,p)(-1)^{k}\frac{1}{k!}\big(e_{\lambda}(t)-1\big)^{k}=\sum_{k=p}^{\infty}L(k,p)(-1)^{k}\sum_{n=k}^{\infty}S_{2,\lambda}(n,k)\frac{t^{n}}{n!}\nonumber\\
&=\sum_{n=p}^{\infty}\bigg(\sum_{k=p}^{n}L(n,k)(-1)^{k}S_{2,\lambda}(n,k)\bigg)\frac{t^{n}}{n!}. \nonumber
\end{align}
Therefore, by \eqref{54} and \eqref{55}, we obtain the following theorem.
\begin{theorem}
For $n,p\in\mathbb{Z}$ with $n\ge p\ge 0$, we have
\begin{displaymath}
	(-1)^{n}S_{2,-\lambda}(n,p)=\sum_{k=p}^{n}L(n,k)(-1)^{k}S_{2,\lambda}(n,k).
\end{displaymath}	
\end{theorem}
Note that
\begin{align}
&\frac{1}{p!}\bigg(\frac{1-e_{\lambda}(t)}{e_{\lambda}(t)}\bigg)^{p}=(-1)^{p}e_{-\lambda}^{p}(-t)\frac{1}{p!}\Big(e_{\lambda}(t)-1\Big)^{p} \label{56} \\
&=(-1)^{p}e_{-\lambda}^{p}(-t)\sum_{k=p}^{\infty}S_{2,\lambda}(k,p)\frac{t^{k}}{k!}=(-1)^{p}\sum_{l=0}^{\infty}(p)_{l,-\lambda}(-1)^{l}\frac{t^{l}}{l!}\sum_{k=p}^{\infty}S_{2,\lambda}(k,p)\frac{t^{k}}{k!}\nonumber \\
&=(-1)^{p}\sum_{n=p}^{\infty}\bigg(\sum_{k=p}^{n}S_{2,\lambda}(k,p)(-1)^{n-k}(p)_{n-k,-\lambda}\binom{n}{k}\bigg)\frac{t^{n}}{n!}.\nonumber
\end{align}
By \eqref{55} and \eqref{56}, we get
\begin{equation}
\sum_{k=p}^{n}L(k,p)(-1)^{k}S_{2,\lambda}(n,k)=(-1)^{p}\sum_{k=p}^{n}S_{2,\lambda}(k,p)(-1)^{n-k}(p)_{n-k,-\lambda}\binom{n}{k}.\label{57}	
\end{equation}
Therefore, by \eqref{57} and Theorem 13, we obtain the following theorem.
\begin{theorem}
For $n,p\in\mathbb{Z}$ with $n\ge p\ge 0$, we have
\begin{displaymath}
	S_{2,-\lambda}(n,p)=(-1)^{p}\sum_{k=p}^{n}S_{2,\lambda}(k,p)(-1)^{k}(p)_{n-k,-\lambda}\binom{n}{k}.
\end{displaymath}	
\end{theorem}
Now, we observe that
\begin{align}
&\frac{1}{p!}\bigg(\log_{\lambda}\bigg(1+\frac{t}{1-t}\bigg)\bigg)^{p}=\frac{1}{p!}\bigg(\log_{\lambda}\bigg(\frac{1}{1-t}\bigg)\bigg)^{p}=\frac{(-1)^{p}}{p!}\Big(\log_{-\lambda}(1-t)\Big)^{p}\label{58} \\
&=(-1)^{p}\sum_{n=p}^{\infty}S_{1,-\lambda}(n,p)(-1)^{n}\frac{t^{n}}{n!}=\sum_{n=p}^{\infty}{n \brack p}_{-\lambda}\frac{t^{n}}{n!}.\nonumber	
\end{align}
On the other hand, by \eqref{4} and \eqref{13}, we get
\begin{align}
&\frac{1}{p!}\Big(\log_{\lambda}\Big(1+\frac{t}{1-t}\Big)\Big)^{p}=\sum_{k=p}^{\infty}S_{1,\lambda}(k,p)\frac{1}{k!}\Big(\frac{t}{1-t}\Big)^{k}\label{59}\\
&=\sum_{k=p}^{\infty}S_{1,\lambda}(k,p)\sum_{n=k}^{\infty}L(n,k)\frac{t^{n}}{n!}=\sum_{n=p}^{\infty}\bigg(\sum_{k=p}^{n}S_{1,\lambda}(k,p)L(n,k)\bigg)\frac{t^{n}}{n!}.\nonumber	
\end{align}
Therefore, by \eqref{58} and \eqref{59}, we obtain the following theorem.
\begin{theorem}
For $n,p\in\mathbb{Z}$ with $n\ge p\ge 0$, we have
\begin{displaymath}
{n\brack p}_{-\lambda}=\sum_{k=p}^{n}S_{1,\lambda}(k,p)L(n,k).
\end{displaymath}	
\end{theorem}

\section{Conclusion}
Carlitz initiated the exploration of degenerate Bernoulli and Euler polynomials, which are degenerate
versions of the ordinary Bernoulli and Euler polynomials.
Along the same line as Carlitz's pioneering work, intensive studies have been done for degenerate versions of quite a few special polynomials and numbers by employing such tools as combinatorial methods, generating functions, umbral calculus techniques, $p$-adic analysis, differential equations, special functions, probability theory and analytic number theory.\par
In this paper, by using generating functions we further studied some identities and properties on the degenerate Stirling numbers of the second kind, the degenerate $r$-Stirling numbers of the second kind, the degenerate Stirling numbers of the first kind and the unsigned degenerate Stirling numbers of the first kind. \par
It is one of our future projects to continue to explore various degenerate versions of many special polynomials and numbers by using aforementioned tools.

\end{document}